\documentclass{amsart}

\usepackage{amsmath,amssymb,amsthm}
\usepackage{graphicx,tikz}
\newtheorem{theorem}{Theorem}

\theoremstyle{definition}

\theoremstyle{remark}

\begin{document}

\title[]{Randomized Kaczmarz converges\\ along small singular vectors}
\thanks{S.S. is supported by the NSF (DMS-1763179) and the Alfred P. Sloan Foundation.}

\author[]{Stefan Steinerberger}
\address{Department of Mathematics, University of Washington, Seattle}
\email{steinerb@uw.edu}

\begin{abstract} Randomized Kaczmarz is a simple iterative method for finding solutions of linear systems $Ax = b$. We point out that the arising sequence $(x_k)_{k=1}^{\infty}$ tends to converge to the solution $x$ in an interesting way: generically, as $k \rightarrow \infty$, $x_k - x$ tends to the singular vector of $A$ corresponding to the smallest singular value. This has interesting consequences: in particular, the error analysis of Strohmer \& Vershynin is optimal. It also quantifies the `pre-convergence' phenomenon where the method initially seems to converge faster. This fact also allows for a fast computation of vectors $x$ for which the Rayleigh quotient $\|Ax\|/\|x\|$ is small: solve $Ax = 0$ via Randomized Kaczmarz. 
\end{abstract}

\maketitle

\section{Introduction}
Let $A \in \mathbb{R}^{m \times n}$ with $m \geq n$ and full rank. Suppose $x \in \mathbb{R}^n$ and  
$$ Ax =b,$$
how would one go about finding $x$?
A popular method for solving these types of problems is the Kaczmarz method \cite{kac} (also known as the Algebraic Reconstruction Technique (ART) in computer tomography \cite{gordon, her, natterer} or as the Projection onto Convex Sets method \cite{cenker, sezan}). Denoting the rows of $A$ by $a_1, \dots, a_m$, the linear system can be written as a collection of inner products
$$ \forall~1 \leq i \leq m: \left\langle a_i, x\right\rangle = b_i.$$
The Kaczmarz method is an iterative scheme: given $x_k$, let us pick an equation, say the $i-$th equation, and modify $x_k$ the smallest possible amount necessary to make it correct: set $x_{k+1} = x_k + \delta a_i,$ where $\delta$ is such that $\left\langle a_i,x_{k+1} \right\rangle = b_i$. Formally,
$$ x_{k+1} = x_k + \frac{b_i - \left\langle a_i, x_k\right\rangle}{\|a_i\|^2}a_i.$$
The Kaczmarz method is then given by cycling through the indices. The method has a geometric interpretation (see Fig. 1): $x_{k+1}$ is obtained by projecting $x_k$ onto the hyperplane given by $\left\langle a_i, x \right\rangle = b_i$. This interpretation is already given in the original 1937 paper of Kaczmarz who ends his paper by saying that `The convergence of this method is completely obvious from a geometric perspective' \cite{kac}. However, given the intricate underlying geometry, convergence rates are difficult to obtain. In a seminal paper, Strohmer \& Vershynin \cite{strohmer} showed that by choosing the equations in a random order (where the $i-$th equation is chosen with a likelihood proportional to $\|a_i\|^2$), the method converges exponentially: more precisely, they show that
$$ \mathbb{E} \left\| x_k - x \right\|_2^2 \leq \left(1 - \frac{1}{\|A\|_F^2 \| A^{-1}\|_2^2}\right)^k \|x_0 - x\|_2^2,$$
where $\|A^{-1}\|_2$ is the operator norm of the inverse and $\|A\|_F$ is the Frobenius norm.
This is a remarkable result in many ways. The analysis of Kaczmarz and Kaczmarz-type algorithms has become quite popular: we refer to 
\cite{bai, eldar, gower, lee, leventhal, liu, ma, need, need2, need3, need4, tan, zouz} and references therein.

\begin{center}
\begin{figure}[h!]
\begin{tikzpicture}[scale=2]
\draw[thick] (-1, -1) -- (1, 1);
\draw[thick] (0, -1) -- (0, 1);
\draw[thick] (1, -1) -- (-1, 1);
\filldraw (0,0) circle (0.03cm);
\node at (0.2, 0) {$x$};
\filldraw (0.7, 0.7) circle (0.03cm);
\node at (0.9, 0.7) {$x_k$};
\filldraw (0, 0.7) circle (0.03cm);
\node at (0.25, 0.78) {$x_{k+1}$};
\draw [ <-] (0.1, 0.7) -- (0.7, 0.7);
\filldraw (-0.7/2, 0.7/2) circle (0.03cm);
\filldraw (-0.35, 0.35) circle (0.03cm);
\node at (-0.6, 0.35) {$x_{k+2}$};
\draw [->] (0, 0.7) -- (-0.3, 0.4);
\end{tikzpicture}
\caption{A geometric interpretation: iterative projections onto the hyperplanes given by $\left\langle a_i, x\right\rangle = b_i$.}
\end{figure}
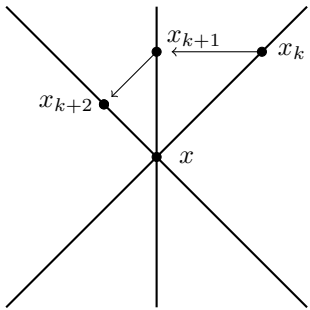
\end{center}

\section{ Results}
\subsection{Summary}
A priori, one might be inclined to believe that the very nature of the algorithm produces an essentially random sequence in the sense that
the normalized vectors $(x_k-x)/\|x_k - x\|$ behave mostly randomly
or, alternatively put, that the approximations converge to $x_k \rightarrow x$ `from all directions'. It is the goal of this short note to show
that this is decidedly \textit{not} the case: the singular vectors $v_1, \dots, v_n$ of the matrix $A$ describe directions of distinguished dynamics. We describe three main observations:
\begin{enumerate}
\item \S 2.1: while $\|x_k-x\|$ decreases exponentially, this is happening at different rates in different subspaces. We show (Theorem 1) that for any (right) singular vector $v_{\ell}$ of $A$ associated to the singular value $\sigma_{\ell}$
$$\mathbb{E}  \left\langle x_{k}-x, v_{\ell}\right\rangle = \left(1 - \frac{\sigma_{\ell}^2 }{\|A\|_F^2}\right)^k\left\langle x_0-x, v_{\ell}\right\rangle.$$
Contributions from the smallest singular vectors ($\sigma_{\ell}$ small) decays the slowest: hence we expect them to be dominant for $k$
large. We recall that for the smallest singular vector $\sigma_n$, we have the identity
$$ \frac{1}{\sigma_n} = \|A^{-1}\|_2 \qquad \mbox{leading to the rate} \qquad 1 - \frac{1}{\|A\|_F^2 \|A^{-1}\|_2^2}$$
which recovers the convergence rate of Strohmer \& Vershynin \cite{strohmer}. Thus in the generic setting one
cannot hope to improve the rate.
\item \S 2.2: this phenomenon also impacts convergence rates:
if $x_k$ is not mainly the linear combination of singular vectors corresponding to small singular values, Randomized Kaczmarz
converges \textit{faster}: we have
$$\mathbb{E} ~ \|x_{k+1} - x\|_2^2 \leq \left(1 - \frac{1}{\|A\|_F^2} \left\| A  \frac{x_k-x}{\|x_k - x\|} \right\|^2 \right) \|x_{k} - x\|_2^2.$$
We observe that this factor may actually be quite large since
$$ \|A^{-1}\|^2\leq  \left\| A  \frac{x_k-x}{\|x_k - x\|} \right\|^2 \leq \|A\|^2.$$
The lower bound is only close to being attained when $x_k-x$ is close to the space spanned by singular vectors given by small singular values of $A$.
\item \S 2.3: once $x_k$ is mainly the combination of singular vectors corresponding to small singular values, the Randomized Kaczmarz Method is also less likely to explore adjacent subspaces. We show that the normalized vector $(x_k-x)/\|x_k-x\|$ then moves only very slowly on $\mathbb{S}^{n-1}$ since
$$ \mathbb{E} ~\left\langle \frac{x_k-x}{\|x_k-x\|}, \frac{x_{k+1}-x}{\|x_{k+1}-x\|}  \right\rangle^2 = 1 - \frac{1}{\|A\|_F^2}\left\| A  \frac{x_k-x}{\|x_k - x\|} \right\|^2.$$
\end{enumerate}

\subsection{Transfer to Small Singular Vectors} We start by establishing the main phenomenon which can be phrased as a simple algebraic identity.
\begin{theorem} Let $v_{\ell}$ be a (right) singular vector of $A$ associated to the singular value $\sigma_{\ell}$. Then
$$\mathbb{E} \left\langle x_{k} - x, v_{\ell} \right\rangle = \left(1 - \frac{\sigma_{\ell}^2 }{\|A\|_F^2}\right)^k\left\langle x_0 - x, v_{\ell}\right\rangle.$$
\end{theorem}
We recall that
$$ \sum_{i=1}^{n}{\sigma_{i}^2} = \|A\|_{F}^2,$$
so the arising factor multiplicative factor is always in between 0 and 1. This shows that the Randomized Kaczmarz method has exponential decay at different
exponential rates that depend on the singular values. In particular, the slowest rate of decay is given by the smallest singular value $\sigma_n$. Since
$$ \frac{1}{\sigma_n} = \|A^{-1}\|_2 \qquad \mbox{the rate is} \qquad 1 - \frac{1}{\|A\|_F^2 \|A^{-1}\|_2^2}$$
which recovers the convergence rate of Strohmer \& Vershynin \cite{strohmer}. In particular, if $x_0$ has contributions from singular vectors corresponding to large singular values, the sequence initially decays faster than expected since components corresponding to larger singular values decay quicker. This has been known for a long time: Randomized Kaczmarz initially acts stronger on large singular values, in particular it serves as a regularizer, we refer to \cite{jin, pic}.
 A previous result in that spirit was obtained by Jiao, Jin \& Lu \cite{jiao} who obtained an `averaged' version of Theorem 1: by decomposing $\mathbb{R}^n = L + H$ into subspaces spanned by low singular vectors $L$ (corresponding to small singular values) and high singular vectors $H$ (corresponding to large singular values), they obtain the estimate
$$ \mathbb{E} \| \pi_{H} x_{k+1} \|^2 \leq (1-c_1) \| \pi_H x_{k}\|^2 + c_2 \| \pi_{L} x_k\|^2$$
where $\pi$ is the projection and $c_1, c_2 > 0$ are given in terms of singular values.

\subsection{Speed of Convergence}
The phenomenon described in Theorem 1 obviously has an effect on the speed of the algorithm: whenever $x_k - x$ is not close to singular vectors of $A$ that correspond to small singular vector, the Randomized Kaczmarz method actually converges faster. Theorem 2 is already implicit in Strohmer \& Vershynin \cite{strohmer} where a formulation of it appears as an intermediate step. 
\begin{theorem}[see also Strohmer \& Vershynin \cite{strohmer}] We have
$$\mathbb{E} ~ \|x_{k+1} - x\|_2^2 \leq \left(1 - \frac{1}{\|A\|_F^2} \left\| A  \frac{x_k-x}{\|x_k - x\|} \right\|^2 \right) \|x_{k} - x\|_2^2.$$
\end{theorem}
Observe that for any vector $\|v\|=1$, we have
$$ \|Av\|^2 \geq \sigma_n^2 \|v\|^2 = \|A^{-1}\|_{2}^{-2} \cdot  \|v\|_2,$$
where $\sigma_n$ is the smallest singular value. But clearly this last inequality is only sharp or close to sharp if $v$ is close to linear combination of singular vectors corresponding to small singular values: otherwise, the quantity is indeed larger.

\subsection{Stability of Small Singular Vectors.}
We conclude with a result of a geometric flavor. A natural object is the sequence of vectors
$$ \frac{x_k - x}{\|x_k - x\|} \in \mathbb{S}^{n-1}.$$
It essentially describes from which direction the sequence $(x_k)_{k=0}^{\infty}$ approaches the solution $x$. The result in the prior sections have shown that we expect this vector will eventually be concentrated on subspaces corresponding to singular vectors associated to small singular values. A different questions would be how much we expect the normalized vector to `move around'.  
\begin{theorem}
If $x_k \neq x$ and $\mathbb{P}(x_{k+1} = x)=0$, then
$$ \mathbb{E} ~\left\langle \frac{x_k-x}{\|x_k-x\|}, \frac{x_{k+1}-x}{\|x_{k+1}-x\|}  \right\rangle^2 = 1 - \frac{1}{\|A\|_F^2}\left\| A  \frac{x_k-x}{\|x_k - x\|} \right\|^2.$$
\end{theorem}
We recall that for two unit vectors $\|u\| = 1 = \|v\|$, the inner quantity $\left\langle u, v\right\rangle^2$ can be interpreted as a notion of distance between the vectors: it is only close to 1 if $u \sim \pm v$. Since we only care about subspaces, the ambiguity in the sign does not play a role. We observe that the direction $(x_k-x)/\|x_k-x\|$ undergoes large fluctuations when $\|A (x_k-x)\|$ is large and is quite stable when $\|A(x_k-x)\|$ is small. Once $x_k$ is mainly comprised of small singular vectors, its direction does not change very much anymore: randomized Kaczmarz spends more time close to these subspaces than it does in other regions of space: the smaller the smallest singular value, the more pronounced phenomenon. Could this property be used for convergence acceleration?

\subsection{Finding Small Rayleigh Quotients}
We discuss an amusing application: finding vectors $\|x\|=1$ such that $\|Ax\|$ is small. In light of our results, such a sequence of vectors can be obtained from the Randomized Kaczmarz method applied to the equation 
$$ Ax = 0.$$
We illustrate this with a simple example: we take $A \in \mathbb{R}^{1000 \times 1000}$ by picking each entry to be iid $\sim \mathcal{N}(0,1)$, we add $100\cdot \mbox{Id}_{1000}$ and then set the last row $a_{1000}$ to be a tiny perturbation of $a_{999}$ by adding 0.01 to each entry. After that, we normalize all rows to have size $\|a_i\|=1$. This generates a matrix having singular 
values between $\sim 0.5$ and $\sim 1.5$ and one tiny singular value at scale $\sim 10^{-4}$. We start with $x_0 = (1,\dots, 1)$ and plot $\|A x_k\|/\|x_k\|$ (see Fig. 2).
We observe that within $5000 \sim 5n$ iterations, the ratio is actually smaller than any singular value but the last one -- this implies that $x_k$ is actually already mostly given by $v_{1000}$. After 10.000 iterations, we observe $|\left\langle x_{10^4}/\|x_{10^4}\|, v_{1000}\right\rangle| \sim 0.9$. 
\begin{center}
\begin{figure}[h!]
\begin{tikzpicture}[scale=2]
\node at (0,0.3) {\includegraphics[width=0.6\textwidth]{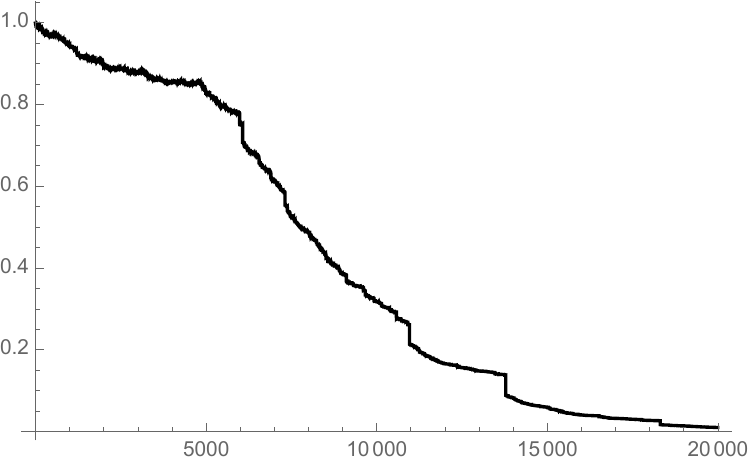}};
\node at (-2.3,1) {\LARGE $\frac{\|A x_k\|}{\|x_k\|}$};
\node at (1.8,-1) {$k$};
\end{tikzpicture}
\caption{A sample evolution of $\|Ax_k\|/\|x_k\|$.}
\end{figure}
\end{center}

\section{Proofs}
\subsection{Reduction to $Ax = 0$} We start with the ansatz
$$ x_k = x + y_k.$$
Plugging in, we obtain
\begin{align*}
x + y_{k+1} &= x_{k+1} = x_k + \frac{ b_i - \left\langle a_i, x_k\right\rangle}{\|a_i\|^2} a_i \\
&=x + y_k + \frac{ b_i - \left\langle a_i, x + y_k\right\rangle}{\|a_i\|^2} a_i \\
&= x + y_k  - \frac{  \left\langle a_i,  y_k\right\rangle}{\|a_i\|^2} a_i + \left( \frac{ b_i - \left\langle a_i, x\right\rangle}{\|a_i\|^2} a_i\right).
\end{align*}
We note that the final term in parenthesis vanishes because $x$ is the exact solution and thus $b_i - \left\langle a_i, x\right\rangle=0$. This results in
$$ y_{k+1} = y_k - \frac{  \left\langle a_i,  y_k\right\rangle}{\|a_i\|^2} a_i$$
which is exactly the Randomized Kaczmarz method for $Ay = 0$. Linearity of all the underlying objects (or, alternatively, the geometric interpretation in terms of hyperplane intersections) explains why this is perhaps not too surprising. In any case, this allows us to restrict ourselves to the case $Ax = 0$.
\subsection{Proof of Theorem 1}
\begin{proof}
As discussed above, we can assume without loss of generality that $b=0$. It remains to understand the behavior of the Randomized Kaczmarz method when applied to $Ax=0$. We have
\begin{align*}
\mathbb{E} \left\langle x_{k+1} , v_{\ell} \right\rangle &=
    \sum_{i=1}^{m}  \frac{\|a_i\|^2}{\|A\|_F^2} \left\langle    x_k - \frac{\left\langle a_i, x_k \right\rangle }{\|a_i\|^2} a_i  , v_{\ell} \right\rangle \\
   &=   \frac{1}{\|A\|_F^2}\sum_{i=1}^{m} \left( \|a_i\|^2 \left\langle x_k, v_{\ell}\right\rangle - \left\langle a_i, x_k\right\rangle \left\langle a_i, v_{\ell} \right\rangle \right)\\
   &= \left\langle x_k, v_{\ell}\right\rangle -   \frac{1}{\|A\|_F^2}\sum_{i=1}^{m} \ \left\langle a_i, x_k\right\rangle \left\langle a_i, v_{\ell} \right\rangle. 
   \end{align*}
This sum has a special structure: $\left\langle a_i, x_k\right\rangle$ is simply the $i-$th entry of $A x_k$ and, likewise, $\left\langle a_i, v_{\ell}\right\rangle$ is simply the $i-$th entry of $A v_{\ell}$. This simplifies the sum
$$ \sum_{i=1}^{m} \ \left\langle a_i, x_k\right\rangle \left\langle a_i, v_{\ell} \right\rangle = \left\langle A x_k, A v_{\ell} \right\rangle.$$
 We then expand $x_k$ into right-singular vectors $v_i$, make use of $A v_i = \sigma_i u_i$ and the orthogonality of the left singular vectors $u_i$ to obtain
   \begin{align*}
 \left\langle x_k, v_{\ell}\right\rangle -   \frac{1}{\|A\|_F^2}\sum_{i=1}^{m} \ \left\langle a_i, x_k\right\rangle \left\langle a_i, v_{\ell} \right\rangle   &= \left\langle x_k, v_{\ell}\right\rangle - \frac{1}{\|A\|_F^2} \left\langle A x_k, A v_{\ell}\right\rangle \\
      &= \left\langle x_k, v_{\ell}\right\rangle- \frac{1}{\|A\|_F^2} \left\langle A \left(\sum_{i=1}^{n}{\left\langle x_k, v_i \right\rangle v_i}\right), A v_{\ell}\right\rangle \\
      &= \left\langle x_k, v_{\ell}\right\rangle- \frac{1}{\|A\|_F^2} \left\langle  \sum_{i=1}^{n}{\sigma_i \left\langle x_k, v_i \right\rangle u_i}, \sigma_{\ell} u_{\ell}\right\rangle\\&= \left\langle x_k, v_{\ell}\right\rangle\left(1 - \frac{\sigma_{\ell}^2 }{\|A\|_F^2}\right)
      \end{align*}
 which is the desired result.
\end{proof}
\subsection{Proof of Theorem 2}
\begin{proof} This proof follows the argument of Strohmer \& Vershynin \cite{strohmer}. If $Z$ is the (vector-valued) random variable assuming the value
$Z = a_j/\|a_j\|^2$ with probability $\|a_j\|^2/\|A\|_F^2$,
then Strohmer \& Vershynin have established the inequality
$$ \|x_{k+1} - x\|_2^2 \leq \left(1 - \left| \left\langle \frac{x_k-x}{\|x_k - x\|}, Z \right\rangle\right|^2 \right) \|x_{k} - x\|_2^2.$$
At this point we simply compute
\begin{align*}
\mathbb{E}  ~\left| \left\langle \frac{x_k-x}{\|x_k - x\|}, Z \right\rangle\right|^2 &= \sum_{i=1}^{m} \frac{\|a_j\|_2^2}{\|A\|_F^2} \left\langle \frac{x_k-x}{\|x_k - x\|},  \frac{a_j}{\|a_j\|_2} \right\rangle^2 \\
&=  \frac{1}{\|A\|_F^2} \sum_{i=1}^{m}   \left\langle \frac{x_k-x}{\|x_k - x\|},  a_j \right\rangle^2  \\
&= \frac{1}{\|A\|_F^2} \left\| A  \frac{x_k-x}{\|x_k - x\|} \right\|^2.
\end{align*}
\end{proof}

\subsection{Proof of Theorem 3}
\begin{proof} As above, we can assume without loss of generality that $b=0$. Moreover, 
since the entire expression is invariant under scaling, we can assume w.l.o.g. that $\|x_k\|=1$. If $\mathbb{P}(x_{k+1}=0)=0$, we have
\begin{align*}
 \mathbb{E} \left\langle x_k, \frac{x_{k+1}}{\|x_{k+1}\|} \right\rangle^2 &= \sum_{i=1}^{m}  \frac{\|a_i\|^2}{\|A\|_F^2} \left\langle x_k,
   \frac{x_k - \frac{\left\langle a_i, x_k \right\rangle }{\|a_i\|^2} a_i}{\| x_k - \frac{ \left\langle a_i, x_k \right\rangle }{\|a_i\|^2} a_i \|}  \right\rangle^2 \\
   &= \sum_{i=1}^{m}  \frac{\|a_i\|^2}{\|A\|_F^2} \frac{ \left\langle x_k,
   x_k - \frac{\left\langle a_i, x_k \right\rangle }{\|a_i\|^2} a_i  \right\rangle^2}{\| x_k - \frac{ \left\langle a_i, x_k \right\rangle}{\|a_i\|^2} a_i \|^2}.
   \end{align*}
   We observe that 
   $$x_k - \frac{ \left\langle a_i, x_k \right\rangle}{\|a_i\|^2} a_i  $$
   is the projection of $x_k$ onto the subspace $\left\langle x, a_i \right\rangle =0$ and thus
   $$ \left\langle x_k, x_k - \frac{ \left\langle a_i, x_k \right\rangle}{\|a_i\|^2} a_i  \right\rangle = \left\|x_k - \frac{ \left\langle a_i, x_k \right\rangle}{\|a_i\|^2} a_i  \right\|^2.$$
   Thus
   \begin{align*}
 \sum_{i=1}^{m}  \frac{\|a_i\|^2}{\|A\|_F^2} \frac{ \left\langle x_k,
   x_k - \frac{\left\langle a_i, x_k \right\rangle }{\|a_i\|^2} a_i  \right\rangle^2}{\| x_k - \frac{ \left\langle a_i, x_k \right\rangle}{\|a_i\|^2} a_i \|^2} 
   &=  \sum_{i=1}^{m}  \frac{\|a_i\|^2}{\|A\|_F^2} \frac{ \|
   x_k - \frac{\left\langle a_i, x_k \right\rangle}{\|a_i\|^2} a_i \|^4}{\| x_k - \frac{\left\langle a_i, x_k \right\rangle }{\|a_i\|^2} a_i \|^2}\\
   &=  \sum_{i=1}^{m}  \frac{\|a_i\|^2}{\|A\|_F^2}  \left\|  x_k - \frac{\left\langle a_i, x_k \right\rangle}{\|a_i\|} \frac{a_i}{\|a_i\|} \right\|^2.
  \end{align*}
  We recall that $\|x_k\| = 1$ and thus
  \begin{align*}
   \left\|  x_k - \frac{\left\langle a_i, x_k \right\rangle}{\|a_i\|} \frac{a_i}{\|a_i\|} \right\|^2 &= \|x_k\|^2 - 2\left\langle \frac{\left\langle a_i, x_k \right\rangle}{\|a_i\|} \frac{a_i}{\|a_i\|}, x_k\right\rangle + \left\| \frac{\left\langle a_i, x_k \right\rangle}{\|a_i\|} \frac{a_i}{\|a_i\|} \right\|^2\\
   &= 1- 2\frac{\left\langle a_i, x_k \right\rangle^2}{\|a_i\|^2} + \frac{\left\langle a_i, x_k \right\rangle^2}{\|a_i\|^2} =  1- \frac{\left\langle a_i, x_k \right\rangle^2}{\|a_i\|^2}. 
  \end{align*}
Altogether, we have
  \begin{align*}
  \sum_{i=1}^{m}  \frac{\|a_i\|^2}{\|A\|_F^2}  \left\|  x_k - \frac{\left\langle a_i, x_k \right\rangle}{\|a_i\|} \frac{a_i}{\|a_i\|} \right\|^2  &= \sum_{i=1}^{m}  \frac{\|a_i\|^2}{\|A\|_F^2}  \left( 1 - \frac{\left\langle a_i, x_k \right\rangle^2}{\|a_i\|^2} \right) \\
   &= \frac{1}{\|A\|_F^2} \sum_{i=1}^{m} \left(\|a_i\|^2 - \left\langle a_i, x_k \right\rangle^2\right) \\
   &= 1 - \frac{1}{\|A\|_F^2} \sum_{i=1}^{m}  \left\langle a_i, x_k \right\rangle^2 = 1 - \frac{\|Ax_k\|^2}{\|A\|_F^2}.
    \end{align*}
\end{proof}

\end{document}